\documentclass[preprint]{imsart}

\RequirePackage[OT1]{fontenc}
\RequirePackage{amsthm,amsmath,booktabs}
\RequirePackage[colorlinks,citecolor=blue,urlcolor=blue]{hyperref}
\PassOptionsToPackage{numbers, compress}{natbib}

\usepackage{amsthm,amsmath,natbib,amssymb,enumerate}
\usepackage{times}
\usepackage{amsthm,amsmath,natbib,amssymb,enumerate}
\usepackage{times}
\usepackage{graphicx} 
\usepackage{subcaption} 

\usepackage{natbib}

\usepackage{algorithm}
\usepackage{algorithmic}

\newcommand{\minimize}{\mathop{\mathrm{minimize}}}

\newcommand{\cA}{\mathcal{A}}
\newcommand{\cB}{\mathcal{B}}

\newcommand{\cS}{\mathcal{S}}

\newcommand{\Ee}{\mathbb{E}}
\newcommand{\cE}{\mathcal{E}}
\newcommand{\Pp}{\mathbb{P}}

\newcommand{\hz}{\hat{z}}
\newcommand{\hE}{\hat{E}}
\newcommand{\hbeta}{\hat{\beta}}

\newcommand{\hB}{\hat B}

\newcommand{\solutionset}{{\cal S}}

\newcommand*\diff{\mathop{}\!\mathrm{d}}

\providecommand{\argmin}{\mathop\mathrm{arg min}}

\providecommand{\sign}{\mathop\mathrm{sign}}

\def\supp#1{\mathrm{supp}({#1})}
\def\unif{\textnormal{Unif}} 

\newcommand{\real}{\mathbb{R}}

\newtheorem{theorem}{Theorem}
\newtheorem{corollary}[theorem]{Corollary}
\newtheorem{lemma}[theorem]{Lemma}

\begin{document} 
\begin{frontmatter}

\title{MAGIC: a general, powerful and tractable method for selective inference}
\runtitle{MAGIC}

\begin{aug}
\author{\fnms{Xiaoying} \snm{Tian}\corref{}\ead[label=e1]{xtian@stanford.edu}},  
\author{\fnms{Nan} \snm{Bi}\corref{}\ead[label=e2]{nbi@stanford.edu}}
\and  
\author{\fnms{Jonathan} \snm{Taylor}\ead[label=e3]{jonathan.taylor@stanford.edu}\thanksref{t1}
}
\runauthor{Tian et al.}

\affiliation{$^1$Stanford University}

\address{Department of Statistics\\ Stanford University 
\\ Sequoia Hall \\ Stanford, CA 94305, USA \\ \printead{e1} \\
\printead*{e2} \\ \printead*{e3}} 

\thankstext{t1}{Supported in part by NSF grant DMS 1208857 and
AFOSR grant 113039.}
\end{aug}

\begin{abstract} 
    Selective inference is a recent research topic that tries to perform
    valid inference after using the data to select a reasonable statistical
    model.  
    We propose MAGIC, a new method for selective inference that is general,
    powerful and tractable. MAGIC is a method for selective inference after
    solving a convex optimization problem with smooth loss and $\ell_1$ penalty. 
    Randomization is incorporated into the optimization problem to boost statistical power. 
    Through reparametrization, MAGIC reduces the problem into a sampling problem
    with simple constraints. MAGIC applies to many $\ell_1$ penalized optimization
    problem including the Lasso,
    logistic Lasso and neighborhood selection in graphical models, all of which 
    we consider in this paper.
\end{abstract} 



\end{frontmatter}

\section{Introduction}
\label{sec:introduction}
There are a great deal of sophisticated statistical learning methods that allow us to 
search through a large number of models and look for meaningful patterns. Having done this search,
we naturally want to judge the apparent associations that have been found.
This has spawned a new area of research called selective inference
\citep{optimal_inference,exact_lasso,selective_sqrt,randomized_response}. Loosely speaking selective inference
recognizes the inherent selection biases in reporting the most ``significant'' results from various statistical models and 
attempts to adjust for the bias. 

At a high level, selective inference involves two stages: First, solve a convex
optimization problem, usually some penalized loss function. Second, perform 
inference in the statistical model suggested by the result of the optimization
problem. For example, we first use the data to solve the Lasso problem, and
then want to form confidence intervals for the variables that are nonzero in the Lasso solution. 
Adjustment for selection results in some constraints on the underlying
distribution. Although various such problems have been studied, most of the papers
only focus on one specific optimization problem. This is necessary as different loss
functions in the optimization problems result in different geometry of the constraints. In this paper,
we introduce a method, called "MAGIC", Monte-carlo Algorithm for General Inference
with Constraints, which provides valid selective inference for optimization 
problems with any smooth loss functions. The advantage of MAGIC compared to 
previous selective inference methods are generality, statistical power and tractability. 
We elaborate each in the following passage:

{\bf Generality:} The generality of MAGIC lies in two aspects: arbitrary
smooth loss function in the penalized optimization problem and the data
distribution from any exponential family. In comparison, the authors in \cite{exact_lasso} considered only 
inference after solving Lasso; The work \cite{optimal_inference} considered 
some exponential families with simple selection rules, but also noticed
the difficulty for inference after solving more complex optimization problems. Finally, 
the work \cite{likelihood_inference} is the closest in generality to this work,
but shows substantially weaker statistical power, which we discuss below.

{\bf Statistical power:} Earlier work \citep{exact_lasso, likelihood_inference} has
provided valid inference after selection, but sometimes lacks power.
Other work
\citep{randomized_response, optimal_inference} suggested introducing randomness
in the optimization algorithm which produces much improved power. 
This is the approach we take in this work. In simulation, we show that MAGIC
produces much improved power over \cite{exact_lasso, likelihood_inference}.

{\bf Tractability:} The earlier work \citep{exact_lasso, likelihood_inference} computes valid
p-values in closed form, thus involving the least computation cost. The framework in 
\cite{optimal_inference, randomized_response} involves sampling in a constrained 
subset in the sample space. Both work used hit-and-run algorithm proposed in 
\cite{hit-and-run}, which is a method to generate distributions on a subset 
of the space. The constrained subsets as described in \cite{optimal_inference,
randomized_response, exact_lasso, likelihood_inference} can be quite complicated and 
depend on the particular loss function.
Algorithms that do not use MCMC, such as \cite{exact_lasso, likelihood_inference},
do not suffer from this problem too much
as they only need to compute the boundary once, but the methods in \cite{optimal_inference,
randomized_response} need to compute the boundary at each step of simulation,
resulting in much more computation cost. MAGIC, however, transforms the 
constrained subset to a canonical set through reparametrization, removing
the computational cost involved in computing the boundary at each step of sampling. 
Thus it is more tractable than \cite{optimal_inference, randomized_response}.

In Section \ref{sec:setup}, we introduce the general form of our randomized
optimization problem, and describe the inference method as well as theory for
MAGIC. Section \ref{sec:example} gives applications of MAGIC to different
statistical learning problems. To demonstrate the applicability of MAGIC,
we give three distinct examples: Lasso, $\ell_1$ penalized logistic regression
and neighbourhood selection \cite{neighbourhood_selection}, which are applied
in regression, classification and Gaussian graphical models respectively.
Section \ref{sec:simulation} includes the comparisons of MAGIC with
existing selective inference methods both in terms of statistical power and
confidence intervals. All proofs are collected in Section \ref{sec:proof} and
the sampling methods are covered in Section \ref{sec:sampling}. We conclude with
discussions about future work in Section \ref{sec:discussion} 

\subsection{Related works}
Most of the theoretical work on high-dimensional data focuses on 
consistency, either the consistency of solutions \citep{negahban_unified_2010,
van2008high} or the consistency of the models \citep{wainwright2009sharp, 
lasso_consistency}.

In the post selection literature, the authors in \cite{posi} proposed the PoSI approach,
which reduce the problem to a simultaneous inference problem. Because of
the simultaneity, it prevents data snooping from any selection procedure,
but also results in more conservative inference. In addition, the PoSI
method has extremely high computational cost, and is only applicable when
the dimension $p < 30$ or for very sparse models. 
The authors \citep{multi_split} proposed a method for computing
p-values that controls false discovery rate (FDR) among all variables. 
This is quite different from the hypothesis testing framework of this work,
as the hypotheses tested in selective inference are chosen as a function
of the data. Hence, the hypotheses tested
are not directly comparable. Furthermore, compared with \cite{multi_split},
MAGIC has the advantage of being able to construct confidence intervals for
the selected variables.

\section{Randomized selective inference}
\label{sec:setup}
\subsection{A randomized selection algorithm}

Many statistical learning problems can be cast as convex optimization problems.
Specifically, we solve the following randomized convex optimization.
\begin{equation}
\label{eq:canonical:random:program}
\hbeta = \minimize_{\beta \in \real^p} \ell(\beta;S) + \lambda \|\beta\|_1 - \omega^T\beta,
\end{equation}
where data $S \sim F$, 
$\ell$ can be the negative log-likelihood for $F$, but generally just needs to be
some convex loss function in $\beta$, the randomization variable
$\omega \sim G$, a distribution on $\real^p$ independent of $F$, $\lambda$,
is fixed. Without randomization, that is to set $G = \delta_0$, the
point mass at $0$, \eqref{eq:canonical:random:program} includes many classical
statistical learning problems, e.g. lasso \citep{lasso}, elastic net
\citep{elastic_net}, $\ell_1$ penalized logistic regression, neighbourhood
selection \citep{neighbourhood_selection}. Although it might seem strange to add
noise to data for model selection, it is seen in other forms in literature and
applications. Common use of data splitting is an example \citep{data_splitting,
wasserman_data_splitting}, as a random subset of data is used for model
selection.  The form of our randomization is also related to
\citep{stability_selection}.  We can control the amount of randomization through the
variance of $G$, usually just a little randomization will produce much improved power.

We define the {\em variable selection map} as
$$
\hat{E}(S, \omega) = \supp{\hbeta(S, \omega)}.
$$ 
For the observed data $S_{obs}$ and an instance of $\omega_{obs}$ both
considered fixed, we define $E = \hat{E}(S_{obs}, \omega_{obs})$ which
is the active set of \eqref{eq:canonical:random:program} and consider
it fixed hereafter. 

After having solved the above problem, we now consider inference for parameters
chosen on the basis of this set of non-zero coefficients $E$. 
Suppose the data $S \sim F$ is a member of an exponential family with parameters $b \in \real^p$
and sufficient statistics $T(S) \in \real^p$. In particular, its density $f_{b}(s)$ has 
the following form,
$$
\frac{df_{b}}{d\mu}(s) = \exp(b^T T(s) - A(b)) 
$$
where $\mu$ is the reference measure on the sample space of $S$ 
and $A$ is the normalizing constant  with $\mu, A$ known.
Having observed a set of selected variables $E$, we can and often do then consider a
submodel of the above model with $b_{-E}=0$. If $E \supseteq \supp{b}$, then our
model is correctly specified. This is the scenario we always consider hereafter. 
For treatment of misspecified models, see \citep{exact_lasso, optimal_inference}. 
Under this submodel, the joint distribution of $(S,\omega)$ is fully specified.
Our target of inference is now $b_E$.

Since $E$ is not given a priori, but selected by the data, it seems to be only fair to
consider $(S,\omega)$ such that $\hat{E}(S,\omega) = E$.
This is equivalent to condition on the event $\{(S, \omega):\hat{E}(S, \omega) = E\}$. 
This is the general approach taken in \citep{exact_lasso, optimal_inference, randomized_response}
to provide valid (selective) inference in the above model.

Let 
$\cA$ be the region where
$\{(S, \omega) \in \cA\} \iff \{\hat{E}(S, \omega) = E\}$,
then this general approach to selective inference requires us to describe the conditional distribution
\begin{equation}
\label{eq:random:inference}
S \mid (S, \omega) \in \cA, \qquad (S, \omega) \sim F \times G.
\end{equation}

We first state the following result,
\begin{theorem} 
\label{thm:pvalue}
Suppose $(S, \omega) \sim F \times G$, $F$ is the exponential family specified above,
with the parameters $b$ satisfying $\supp{b} \subseteq E$. 
$G$ is a distribution on $\real^p$ and $\cA$ is defined as above.
Then for any variable $j \in E$, there exists a p-value function $P_j: \supp{F} \rightarrow [0,1]$, such that
\begin{equation}
\label{eq:alpha_level}
\Pp_{F \times G} \left[P_j(S; \cA) \leq \alpha \mid (S, \omega) \in \cA \right] \leq \alpha,
\end{equation}
under the null hypothesis $H_{0j}: b_j = 0$. The function $P_j$ only depends on data $S$ and $\cA$. 
\end{theorem} 

In some cases, equality holds in \eqref{eq:alpha_level}, we will discuss the 
conditions in the proof. In this case, the test proposed above is the 
Uniformly Most Powerful Unbiased test \citep{optimal_inference}, providing 
theoretical ground for the power of MAGIC.
Theorem \ref{thm:pvalue} gives a construction of the p-value, which we can
use to reject the null hypothesis at level $\alpha$.
We will give the exact construction of $P_j$ in the proof of Theorem \ref{thm:pvalue}, which is an
multivariate integral and is hard to compute in general. We instead try to acquire
samples from \eqref{eq:random:inference} and approximate the multivariate
integral. The constrained region $\cA$ is the bottleneck for the sampling,
which is complicated and specific to the loss function $\ell$. However,
through a reparametrization, we can form the constrained region as a simple
set that is independent of $\ell$.

\subsection{Augmented parameter space}
Once we solve the optimization \eqref{eq:canonical:random:program}, 
we get $\hbeta$ the solution and $\hz$ the subgradient of $\|\hbeta\|_1$.
$\hbeta$, $\hz$ are functions of $(S, \omega)$. 
We formally define the {\em optimization map} as follows:
\begin{equation}
\label{eq:opt:map}
(s, \omega) \overset{\hat{\theta}}{\mapsto} (s, \hat{\beta}(s, \omega), \hat{z}(s, \omega)) \in \solutionset^F(\ell),
\end{equation}
where
\begin{equation}
\label{eq:canonical:set}
\solutionset^F(\ell) \overset{\text{def}}{=} \biggl \{(s, \beta, z): 
  s \in \text{supp}(F), ~ 
  \ell(\beta;s) < \infty,~ 
  \|\beta\|_1 < \infty, ~
  z \in \partial \|\beta\|_1 \biggr\}.
\end{equation}
$\solutionset^F(\ell)$ is the set of possible values $(s, \beta, z)$ where 
there will be a solution to \eqref{eq:canonical:random:program}. We call $\solutionset^F(\ell)$
the augmented parameter space. Note $\hbeta$ and $\hz$ are random variables (through the randomness
in $(S, \omega)$). One way to describe their distribution,
is to find the inverse of the map $\hat{\theta}$ to reconstruct $\omega$ from $(S, \hbeta, \hz)$. 
In the following passage, we denote $\hbeta$, $\hz$ as the random variables and $\beta$, $z$ as the
corresponding integration variables when writing out the density.


\subsection{Reconstruction and description of the constrained set}
Let $\hbeta_E$, $\hz_E$ be $\hbeta$ and $\hz$ restricted to $E$, and $\hz_{-E}$ the subgradients restricted to $E^c$. 
To make the notation easier, we define the gradient map $\gamma$ 
$$
\gamma(s, \beta) = \partial_{\beta} \ell(\beta; s).
$$ 
\begin{lemma}
\label{lem:KKT}
Through the reparametrization in the optimization map \eqref{eq:opt:map}, the 
selection event $\{\hat{E}(S,\omega) = E\}$ is equivalent to
\begin{equation}
    \label{eq:KKT}
    \begin{cases}
    \gamma(S, \hbeta)+ \lambda \cdot \hz - \omega = 0, \\
    \text{sign}(\hbeta_E) = \hz_E, \quad \|\hz_{-E}\|_{\infty} \leq 1.
    \end{cases}
\end{equation}
\end{lemma}

Lemma \ref{lem:KKT} provides a reconstruction of $\omega$ using $S$, $\hbeta$
and $\hz$. The {\em reconstruction map} is defined as,
$$
\psi(s,\beta,z) \overset{\text{def}}{=} (s, \gamma(s, \beta) + \lambda \cdot z) = (s, \omega).
$$
It is thus easy to see that the distribution of $(S, \hbeta, \hz)$ follows satisfies
the following distributional condition,
$$
(S, \gamma(S, \hbeta) + \lambda \cdot \hz) \sim F \times G.
$$
 
Moreover, using Lemma \ref{lem:KKT}, the distribution for inference \eqref{eq:random:inference}
can be rewritten as 
\begin{align*}
    S &\mid (\hbeta(S, \omega), \hz(S, \omega)) \in {\cal B}, 
    \qquad (S, \omega) \sim F \times G, \\
{\cal B} &= \bigg\{\hbeta_{-E} = 0, 
    \quad \sign(\hbeta_E) = \hz_E, \quad \|\hz_{-E}\|_{\infty} < 1\bigg\}.
\end{align*}
Note that ${\cal B}$ is a much nicer set than $\cA$ in that it only requires  
$\hbeta_E$ to be in a certain quadrant and $\|\hz_{-E}\|_{\infty} < 1$.

Combining these two observations above, we have the following theorem.
We denote by $T_{E\backslash j} \in \real^{|E|-1}$
the sufficient statistics $T \in \real^p$ restricted to the set $E - \{j\}$.
\begin{theorem}[Sampling for MAGIC]
\label{thm:change:measure}
Through change of variables \eqref{eq:opt:map}, the law for selective inference \eqref{eq:random:inference}
is equivalent to
\begin{equation}
\label{eq:change:measure}
S | (\hbeta, \hz) \in {\cal B}, \quad  (S, \gamma(S, \hbeta) + \lambda \cdot \hz) \sim F \times G.
\end{equation}
Moreover, suppose $F$, $G$ has densities $f$ and $g$ respectively, the joint distribution of 
$(S,\hbeta,\hz)$ has density proportional to
\begin{equation}
\label{eq:change:measure:density}
f(s) \cdot g(\gamma(s, \beta) + \lambda \cdot z) \cdot J\psi(s,\beta,z) \cdot 1_{\cal B}(\beta,z)
\end{equation}
with the Jacobian denoting the determinant of the derivative of the map $\psi$ with respect
to $(\beta,z)$ on the fiber over $s$.

Furthermore,
assuming the assumptions in Theorem \ref{thm:pvalue}, $P_j(S)$ can
be computed (with approximation) with samples from \eqref{eq:change:measure} and
further conditional on the sufficient statistics $T_{E \backslash j}(S)$.
\end{theorem}

Theorem \ref{thm:change:measure} gives the explicit density of the law \eqref{eq:change:measure} 
up to a constant. In the proof we specify how to use the samples from \eqref{eq:change:measure}
to approximate the p-value function $P_j$.
A natural choice of sampling would be the Metropolis-Hastings
method or perhaps the projected Langevin method \citep{projected_langevin}. To condition on $T_{E\backslash j}(S)$, 
we just need to make sure the proposal does not move $T_{E\backslash j}(S)$ in each step.
Such choice of the proposal is usually natural, for examples see Section \ref{sec:sampling}. 
The boundary constraint is ${\cal B}$, which needs small adjustment from 
the original Metropolis-Hasting method. Detailed description is in Section \ref{sec:sampling}.
After acquiring such samples, we can use them to approximate the p-value function $P_j$ in Theorem \ref{thm:pvalue}.

All the previous work on selective inference also conditions on the
observed signs $\hat{z}_E$.
$$
S \mid (S, \omega) \in \cA, \quad \sign\left(\hbeta_E(S, \omega)\right) = z_{E,obs}
$$ 
where $z_{E,obs} = \sign(\hbeta(S_{obs}, \omega_{obs}))$ is considered fixed.
The work \cite{exact_lasso} explains that any inference valid under this law, would be
valid under \eqref{eq:random:inference}. Note the additional constraint simply
requires $\hbeta_E$ to be in the quadrant specified by $z_{E, obs}$. 
In what  
follows, we also condition on $\hat{z}_E$.

\section{Examples}
\label{sec:example}
\subsection{Randomized Lasso}
Consider linear regression setting where data $y \sim N(Xb, \sigma^2 I)$,  
$X \in \real^{n \times p}$ is fixed, $\sigma^2$ is known.
Instead of solving the original Lasso proposed by \cite{lasso}, we solve the following randomized version of it,
\begin{equation}
    \label{eq:random:lasso}
    \hbeta = \argmin_{\beta \in \real^p} \frac{1}{2}\|y - X\beta\|_2^2 + \lambda \|\beta\|_1 - \omega^T \beta.
\end{equation}
The gradient of the loss $\gamma(y, \hbeta) = -X^T(y - X\hbeta)$ and $\hz$ is the subgradient for $\|\hbeta\|_1$.
The reconstruction map $\psi(y, \beta, z) = \big(y, \lambda \cdot z - X^T(y-X\beta) \big)$. 
Suppose $E$ is the active set of \eqref{eq:random:lasso}, then 
we model the data by $F = N(X_Eb_E, \sigma^2 I)$, $S = y$ and
$b_E$ is the target for inference.

\begin{corollary}[Randomized Lasso sampler]
    \label{cor:lasso:sampler}
    If $E \supseteq \supp{b}$, then conditioning on
    $(E,z_{E,obs})$, the joint distribution of $(y, \hbeta_E, \hz_{-E})$ can be used for inference (for $b_E$). 
    Further, it has density proportional to
    \begin{equation}
        \label{eq:lasso:density}
    \exp\left(-\frac{\|y - X_Eb_E\|_2^2}{2 \sigma^2}\right) \cdot g\left(
     \lambda \begin{pmatrix} 
        z_{E,obs} \\
        z_{-E}
    \end{pmatrix} 
    - X^T(y-X_E\beta_E)\right)
    \end{equation}
    supported on $\sign(\beta_E) = z_{E, obs}$ and $\|z_{-E}\|_{\infty} < 1$.
\end{corollary}
We thus can obtain samples $(y, \hbeta_E, \hz_{-E})$ for any $b_E$ in the null
hypothesis and use Theorem \ref{thm:pvalue} and Theorem \ref{thm:change:measure}
to construct valid p-values. Detailed algorithm is specified in Section \ref{sec:sampling}.

\subsection{L1-penalized logistic regression}
\label{sec:logistic_lasso}
In practice, many statistical learning problems are classification problems, e.g. spam classification, tumor 
analysis, etc. Suppose $x_i \overset{iid}{\sim} F_X$, $x_i \in \real^p$, $y_i|x_i \sim \text{Bernoulli}(x_i^T b)$,
$F_X$ is unknown and $p$ fixed, $S = (X, y)$. 
The logistic loss is 
$$
\ell(\beta) 
= -\frac{1}{\sqrt{n}} \left[\sum_{i=1}^n y_i \log \pi(x_i \beta) + (1-y_i) \log(1-\pi(x_i \beta)) \right],
$$
where $\pi(x) = \exp(x) / (1+\exp(x))$.
The randomized logistic regression solves the following problem
\begin{equation}
    \label{eq:randomized_logistic}
    \hat{\beta} = \text{argmin}_{\beta \in \mathbb{R}^p} \ell(\beta) + \lambda\|\beta\|_1 - \omega^T\beta + \frac{\epsilon}{2}\|\beta\|_2^2
\end{equation}
with $\epsilon > 0$ small and fixed.
The addition of the term with $\epsilon$ is to ensure the existence of the solution to \eqref{eq:randomized_logistic}.
We explicitly express the $\epsilon$ term, but in general it can be absorbed into the loss function.

Suppose $E$ is the active set of \eqref{eq:randomized_logistic}, then $b_E$ is the target of inference.
With slight abuse of notation, we allow $\pi: \real^n \to \real^n, x \mapsto \pi(x)$ to be the function
applied on each coordinate of $x \in \real^n$. With some algebra, we have the reconstruction map for $\omega$ 
$$
\omega= \lambda \cdot \hz - \frac{1}{\sqrt{n}}X^T \left[y - \pi(X\hbeta)\right] + \epsilon \hbeta 
$$ 
To sample $(X, y)$ jointly is not feasible when $F_X$ is unknown.
Two observations help us circumvent
it and even make the sampling more efficient. First, the reconstruction map for $\omega$ only involve the random vector 
$$ 
\begin{aligned}
\nabla \ell(\hbeta_E) &= -\frac{1}{\sqrt{n}}X^T(y - \pi(X_E\hbeta_E)) \\
&\approx -\frac{1}{\sqrt{n}}X^T\left(y - \pi(X_E\bar{\beta}_E) - W(X_E\bar{\beta}_E)X_E(\hbeta_E-\bar{\beta}_E)\right)
\end{aligned}
$$
where $\bar{\beta}_E$ is the MLE for the unpenalized logistic regression with only the variables in 
$E$ and $W(X\beta)=\text{diag}(\pi(X\beta)(1-\pi(X\beta)))$
is the weight matrix. Alternatively, we might take $\bar{\beta}_E$ to be the
one-step estimator in the selected model starting from $\hbeta_E$ \citep{likelihood_inference}.
The gradient $\nabla \ell(\hbeta_E)$ can be reconstructed, up to a Taylor remainder, from $\hbeta$ and 
the random vector
$$
T=
\begin{pmatrix}
\bar{\beta}_E \\
 X_{-E}^T( y - \pi(X_E\bar{\beta}_E)) 
\end{pmatrix}.
$$
Moreover, when $p$ is fixed, pre-selection, our random vector $T$ properly scaled is asymptotically normal 
and when the selected model is correct ($E \supseteq \supp{b}$):
\begin{equation}
\label{eq:asymptotic:normality}
\frac{1}{\sqrt{n}} \bigg[T - \begin{pmatrix}b_E \\ 0 \end{pmatrix} \bigg] \overset{D}{\to} N(0, \Sigma)
\end{equation}
where $\Sigma$ is estimable from the data \citep{likelihood_inference}.
Since asymptotically $T$ is from an exponential family with parameters $b_E$, 
Theorem \ref{thm:pvalue} states the p-value is a function of $T$ only.  
Thus instead of sampling $(X, y)$, we only need to sample the distribution $T$.

\begin{theorem}
    \label{thm:logistic}
    Suppose $E \supseteq \supp{b}$ and conditioning on $(E, z_{E,obs})$, 
    the joint distribution of $(T, \hbeta_E, \hz_{-E})$ can be used for inference. 
    Then the distribution of $(T, \hbeta_E, \hz_{-E})$ asymptotically (with $p$ fixed, $n \rightarrow \infty$) has density
    \begin{equation}
        \begin{aligned}
        \label{eq:asymptotic:density}
        \phi(T) \cdot g\bigg(&\frac{1}{\sqrt{n}}\begin{pmatrix} X_E^TW(X_E\bar{\beta}_E) X_E(\hbeta_E - \bar{\beta}_E) \\ X_{-E}^TW(X_E\bar{\beta}_E) X_E(\hbeta_E - \bar{\beta}_E) -X_{-E}^T(y-\pi(X_E\bar{\beta}_E)) \end{pmatrix} \\
        & + \lambda 
        \begin{pmatrix}
            z_{E,obs} \\
            z_{-E}
        \end{pmatrix} + \epsilon 
        \begin{pmatrix}
            \hbeta_E \\
            0
        \end{pmatrix}
        \bigg)
        \end{aligned}
    \end{equation}
    where $\phi$ is the density for $N((b_E,0),  \Sigma)$.

\end{theorem}

\subsection{Neighborhood Selection}
\label{sec:ns}
Gaussian graphical models have recently become a very popular way to study network structures. In particular, it has often been 
used on many types of genome data (e.g. gene expression, metabolite concentrations etc.) 
Suppose the data we observe is $X \in \real^{n \times p}$, where each row of $X$ is independently
distributed as $N(\mu, \Sigma)$, $\mu \in \real^p, ~ \Sigma \in \real^{p \times p}$.

It is of interest to study the conditional independence structure 
of the variables $\{1,2, \dots, p\}$. The conditional independence structure  
is conveniently represented by an undirectional graph $(\Gamma, \cE)$, where the nodes $\Gamma = \{1,2,\dots,p\}$, and there is an edge between 
$(i,j)$ if and only if $x_i \not\perp x_j$ conditional on all the other variables $\Gamma \backslash \{i, j\}$. Moreover, 
assuming the covariance matrix $\Sigma$ is not singular, we denote the inverse covariance matrix $\Theta = \Sigma^{-1}$, then
$$
x_i \perp x_j | X_{\Gamma \backslash \{i, j\}} \iff \Theta_{ij} = 0.
$$
In many applications of Gaussian graphical models, we assume the sparse edge structure, where we can hope to 
recover the edgeset $\cE$ even when $n < p^2$.
The authors in \cite{neighbourhood_selection} proposed neighborhood selection with the Lasso to achieve this goal. 
The algorithm can be formulated as the following optimization problem, for any node $i$
\begin{equation}
    \label{eq:ns_optimize}
\hbeta^{i, \lambda} = \argmin_{\beta: \beta_i = 0} \left(n^{-1} \|x_i - X\beta\|_2^2 + \lambda \|\beta\|_1\right),
\end{equation}
where $x_i$ is the $i$-th column of $X$, $\lambda$ is chosen according to
Chapter 3 of \citep{neighbourhood_selection} and considered fixed. 
Denote $\hB = (\hbeta^1, \hbeta^2, \dots, \hbeta^p)$, we propose the randomized version of \eqref{eq:ns_optimize},
\begin{equation}
    \label{eq:ns_randomize}
    \hB = \argmin_{B: B_{ii} = 0} \|X - XB\|_F^2 + \lambda \|B\|_1 - \Omega B, 
\end{equation}
where $\Omega = (\omega^1, \dots, \omega^p),~ \omega^i \overset{i.i.d}{\sim} G$. 
Let $E^i = \supp{\hbeta^i}$, and $E = (E^1, \dots, E^p)$. Since $E$ is usually not
symmetric, we instead look at the set 
$$
E^{\lor} = \{(i,j) | E_{ij} = 1 \text{ or } E_{ji} = 1\}.
$$
Our target for inference is $\{\Theta_{ij}, ~(i,j) \in E^{\lor}\}$.
Note \eqref{eq:ns_randomize} is the matrix form of \eqref{eq:canonical:random:program},
and the reconstruction maps are decomposable across the $p$ nodes; Therefore, we have the following corollary,
\begin{corollary}
\label{cor:ns_sampling}
    Suppose $E$ is the active set for \eqref{eq:ns_randomize}, and $\hz_{E}$ is the corresponding signs of $\hB_E$,
    then conditioning on $(E, z_{E, obs})$, the distribution of $(X, \hB_E, \hz_{-E})$ can be used for inference. 
    Furthermore, if we assume $\Theta_{ij} = 0, ~ i \neq j \text{ and } (i,j) \not\in E^{\lor}$, then
    the joint distribution of $(X, \hB_E, \hz_{-E})$ has the following density,
    \begin{equation}
        \label{eq:ns_law}
    \begin{aligned}
    &\exp\left[-\frac{1}{2}\sum_{i=1}^p\Theta_{ii}\|x_i\|^2 + \sum_{(i,j) \in E^{\lor}} \Theta_{ij}x_i^T x_j\right] \\
    &\hspace{20pt} \cdot \prod_{i \in \Gamma} g\left(
    \lambda \begin{pmatrix}
        z_{E^i, obs} \\
        z_{-E^i}^i
    \end{pmatrix} - 
    X_{-i}^T (x_i - X_{E^i} \beta^i_{E^i})
    \right)  
    \cdot \det(X_{E^i}^T X_{E^i}).
    \end{aligned}
    \end{equation}
\end{corollary}

\section{Simulation}
\label{sec:simulation}
Theorem \ref{thm:pvalue} states that our p-values should be valid at level $\alpha$, for
any $\alpha \in [0,1]$, see \eqref{eq:alpha_level}. In fact, all the three examples above
satisfy the condition such that the Type-I error for any level-$\alpha$ test would be equal to
(or asymptotically equal to) $\alpha$. That is equivalent as saying the p-values follow $\unif(0,1)$ distribution. 
To validate Theorem \ref{thm:pvalue} and Theorem \ref{thm:change:measure}, 
we ran the following simulations for each of the examples in Section \ref{sec:example}. Our data
is generated as follows, for Lasso,
$$
y \sim N(Xb, \sigma^2 I), ~X \in \real^{n \times p} \text{, fixed, } \|b\|_0 = s,
$$
where $s \ll p$. The framework works for arbitrary $n$ and $p$. To demonstrate the applicability
of our framework in high dimensions, we set $n=50,~ p=100,~ s=7$.  
For logistic Lasso problem, 
$$
x_i \sim N(0, I), \quad y_i|x_i \sim \text{Bernoulli}(\pi(x_i b)),\quad 
\pi = \frac{\exp(x)}{1 + \exp(x)}, \quad \|b\|_0 = s.
$$
The framework for logistic regression is fixed $p$ and $n \rightarrow \infty$. Thus we take $n=500,~ p=50,~ s=5$.
For both of the examples above, the {\em signal to noise ratio} (snr) is $7$.
Finally, for neighborhood selection, the data matrix is $X \in \real^{n \times p}$, each row of $X$ is i.i.d from $N(0, \Theta^{-1})$.
We take $n=100,~ p=30$, note this is a high-dimensional setting since we have $30 \times 30$ unknown parameters. But 
only $1\%$ of off-diagonal elements of $\Theta$ is non-zero, and the non-zero off-diagonal entries of $\Theta$ are taken to be $\rho = 0.245$
and the diagonal elements are $1$. $\rho = 0.245$ is chosen because any value less than $0.25$ would ensure $\Theta$ is positive definite \citep{neighbourhood_selection}.

\begin{figure}[ht]
\includegraphics[width=\columnwidth]{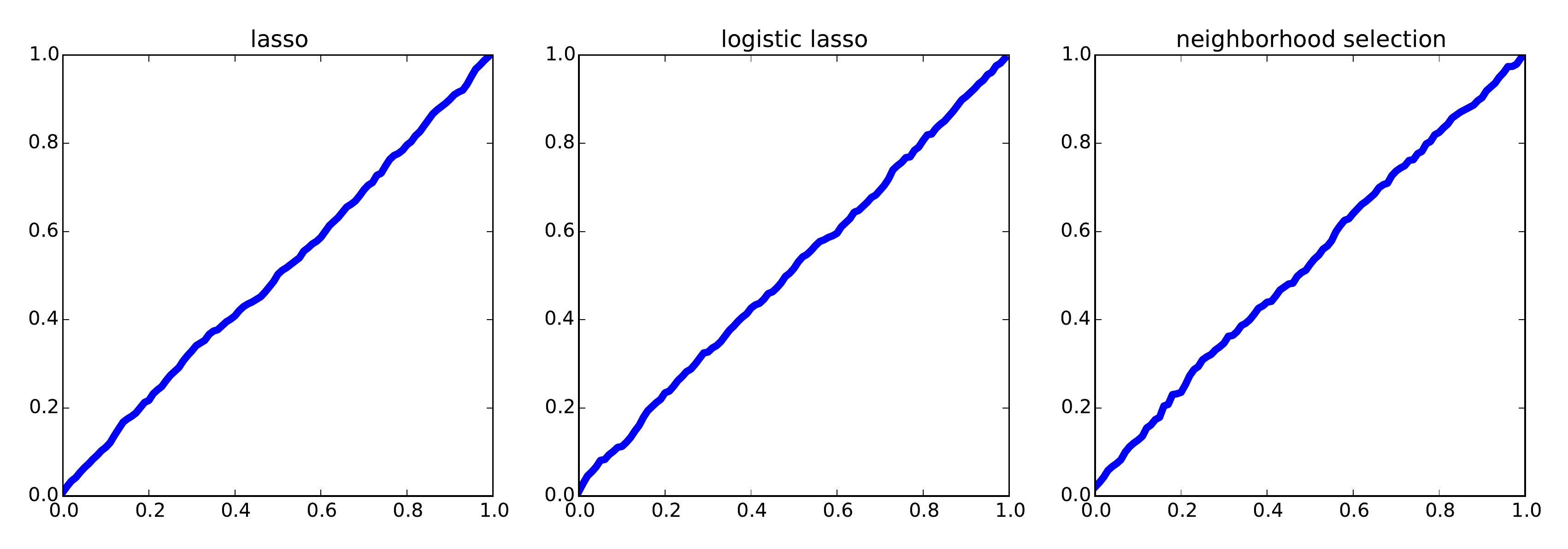}
\caption{Empirical cdf of null p-values, generated from $100$ instances of each problem. 
    We use Laplace noise for randomization.}
\label{fig:pvalues}
\vspace{-2.5mm}
\end{figure}

For each $j \in E$, we test the hypothesis $H_{0j}: b_j = 0$, against a two-sided alternative hypothesis.
We call the p-values the null p-values when the null hypothesis is true and alternative p-values otherwise.
When the active set $E$ (or $E^{\lor}$) from the problem covers $\supp{b}$ (or $\supp{\Theta}$), 
the null p-values should follow $\textrm{Unif}(0,1)$. Figure \ref{fig:pvalues} is the plot for 
the empirical cdf for the null p-values computed from Lasso, logistic Lasso and
neighborhood selection. We see that all the null-pvalues follow the uniform distribution, verifying 
our Theorem \ref{thm:pvalue} and Theorem \ref{thm:change:measure}. 
\subsection{Comparisons of statistical powers}
As we mentioned in Section \ref{sec:introduction}, randomization significantly boosts power.
This is shown in both hypothesis testing and confidence intervals. We describe what it means
in both aspects. For a valid selective level-$\alpha$ test, Type-I error is controlled at $\alpha$
conditional on selection. We hope to achieve valid tests with high power. In the selective
inference framework, statistical power is simply defined as the power in the selected model 
\citep{optimal_inference, randomized_response}. If $E \supseteq \supp{b}$, then for any $j \in E$, 
$$
power = \Pp\big[\text{reject } H_{0j} \mid H_{1j} \text{ is true, } E \text{ is selected}\big]. 
$$
The selective inference framework also offers confidence intervals by inverting a valid test,
for examples, see \cite{exact_lasso}. We want short confidence intervals which
have the desired coverage guarantees. MAGIC enjoys higher statistical power (shorter intervals),
the tradeoff 
is slightly worse selected models as we added randomization for model selection. However, the 
tradeoff is highly in favor of MAGIC. Usually just a small amount of randomization will
dramatically increase statistical power. In the linear regression case, this has been shown in 
\citep{randomized_response} with simulated data. In the following passage, we give numerical 
comparisons on both a real dataset and simulated data. 


\subsubsection{In vitro HIV drug resistance}
\begin{figure}
    \centering
    \caption{Confidence intervals for selected genes in 3TC DATA}\label{fig:3TC}
    \begin{subfigure}[b]{0.478\textwidth}
        \includegraphics[width=\textwidth]{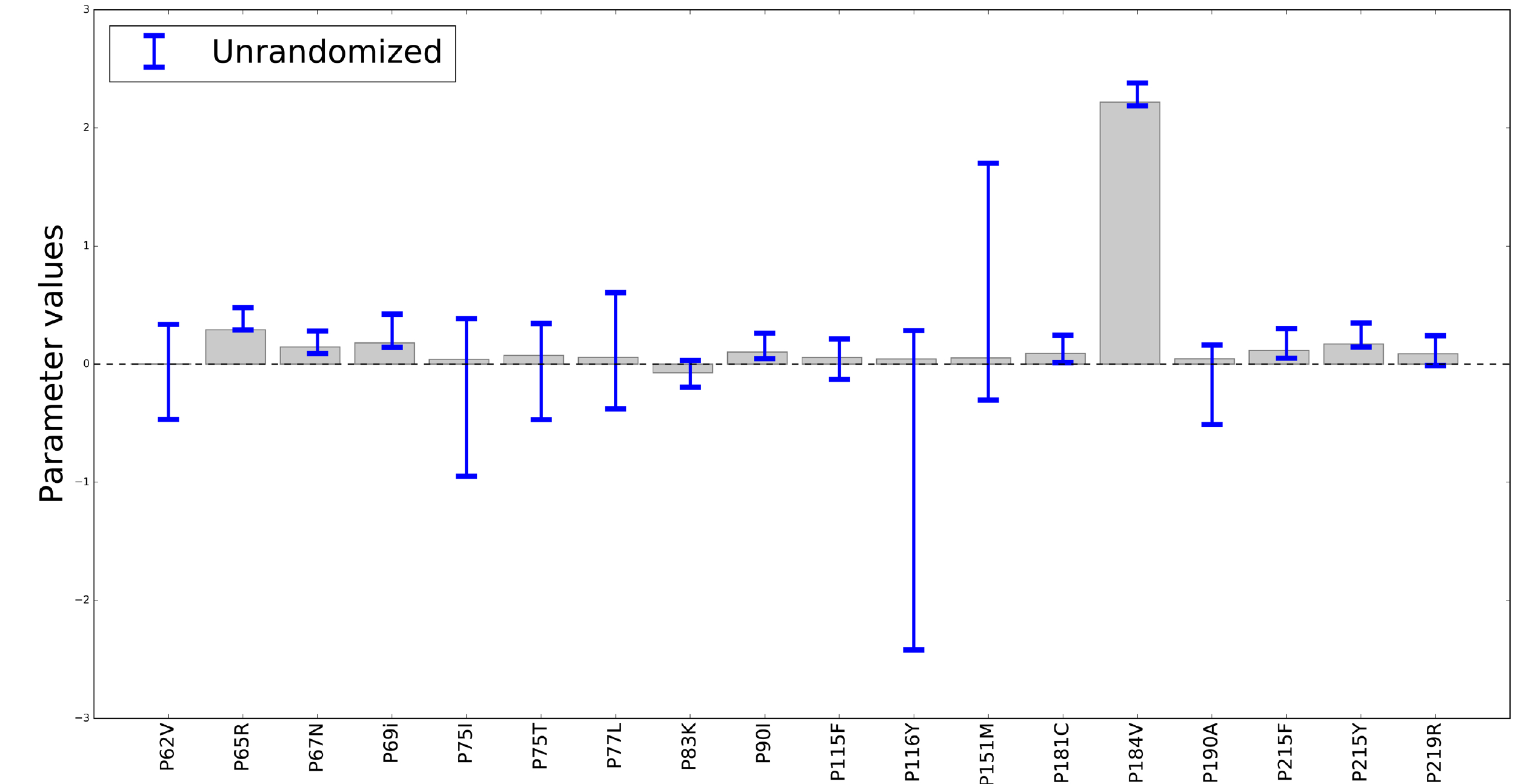}
        \caption{Selective intervals without randomization}
        \label{fig:nonrandomized:3TC}
    \end{subfigure}
    ~ 
    \begin{subfigure}[b]{0.45\textwidth}
        \includegraphics[width=\textwidth]{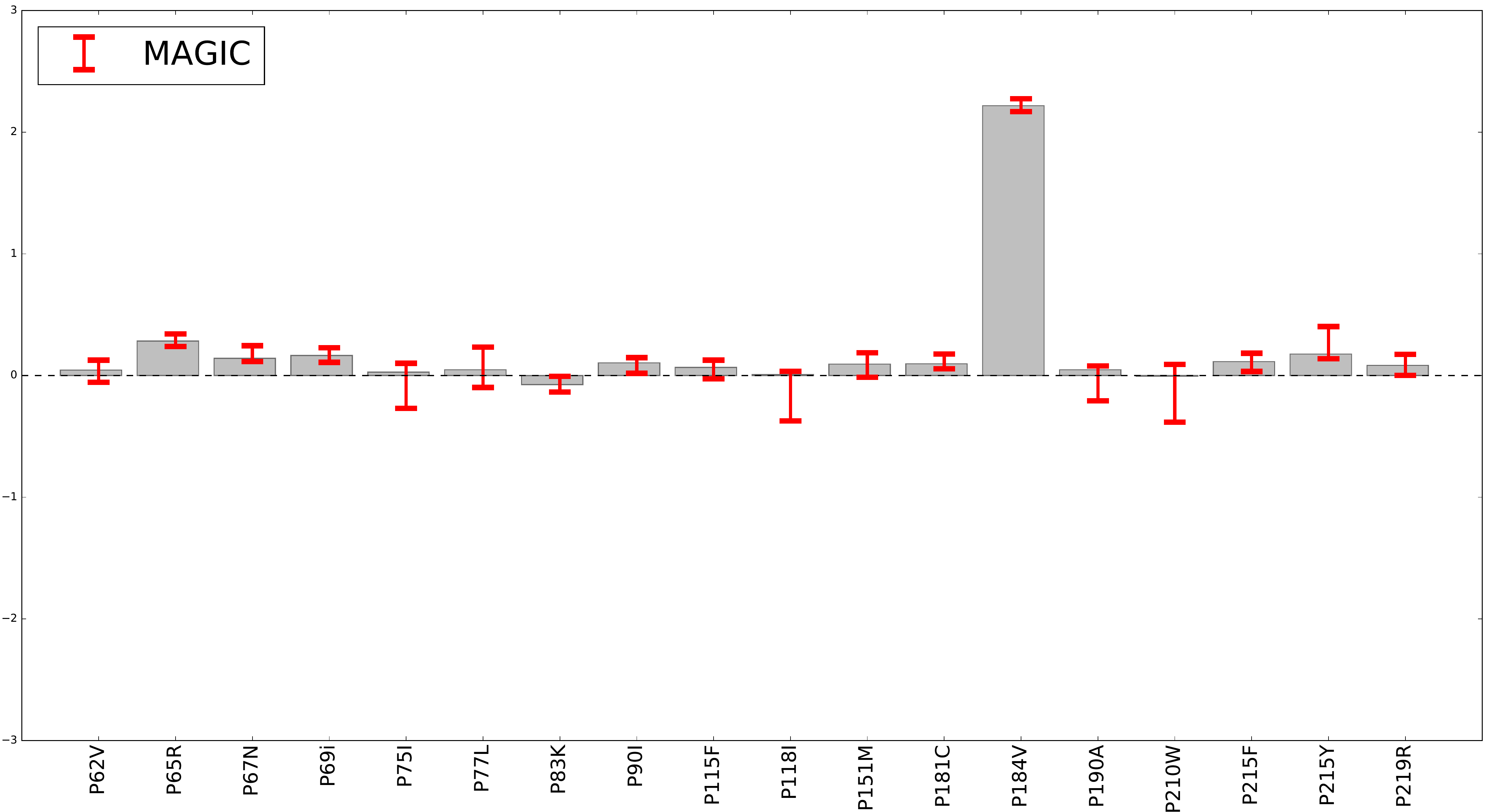}
        \caption{Selective intervals with randomization}
        \label{fig:randomized:3TC}
    \end{subfigure}
\vspace{-5mm}
\end{figure}

In \cite{rhee2006genotypic}, the authors study the genetic basis of drug resistance in HIV, using markers of inhibitor mutations to predict a quantitative measurement of susceptibility to several antiretroviral drugs. 
The hope is to find the mutations highly correlated with the susceptibility to drugs.  
We apply Lasso to the protease inhibitor subset of their data and select the potential mutations set for one of the drugs, Lamivudine (3TC). 
We then compute the OLS estimator in the selected set of gene mutations, and form confidence intervals for the coefficients
(Figure \ref{fig:3TC}). The grey bars are the OLS estimates with only the selected mutations. 
The confidence intervals are adjusted for selection and should have the desired coverage $90\%$.
We report the estimators together with the confidence intervals. The procedure in 
left panel \ref{fig:nonrandomized:3TC} is the same as \cite{exact_lasso} without randomization
in selecting the mutations. The right panel \ref{fig:randomized:3TC} in contrast uses the MAGIC
framework for LASSO with randomization
$\omega \sim N(0, 0.1\sigma_{cv}^2)$, where $\sigma_{cv}$ is the noise level estimated by cross-validation. 
Note the mutations selected by the two methods only differ by $3$ mutations, with small effects, and
the OLS estimator for the common mutations are very close. But the randomized selection
procedure gives much shorter confidence interval across all mutations, demonstrating the 
advantage and practicality of our methods.

\subsubsection{Statistical power comparison with simulated data}
\begin{figure}[ht]
\begin{center}
\caption{The alternative p-values computed from the MAGIC framework highly concentrated 
around $0$, while without randomization the p-values are more evenly distributed between $[0,1]$,
The statistical powers are $0.887$ for MAGIC v.s. $0.176$ for non-randomized
procedure with a level-$0.05$ test.}
\label{fig:pvalue_alt}
\centerline{\includegraphics[width=0.5\columnwidth]{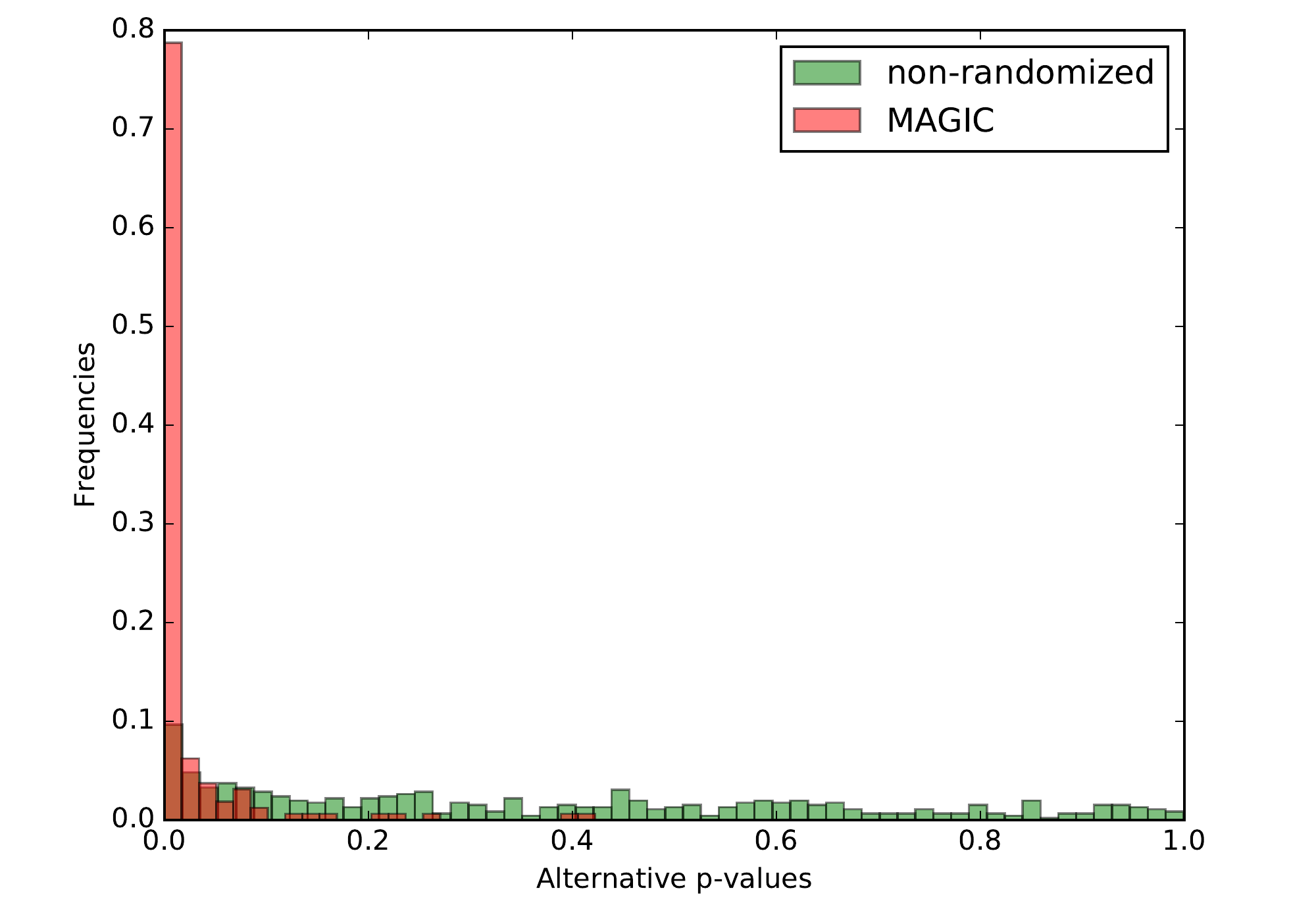}}
\end{center}
\vspace{-5mm}
\end{figure}

In this section, we compare more specifically the tradeoff between power and model
selection using simulated data. The authors in \cite{likelihood_inference}
offered explicit calculations of p-values after the model is selected by 
$\ell_1$ penalized logistic regression or graphical Lasso. Both examples can be
considered in the MAGIC framework. Simulations in \cite{likelihood_inference} showed
that graphical Lasso has worse power than $\ell_1$ penalized logistic regression.
Therefore, we compare our framework to the latter. We assume the same setup
as before, our randomization noise is $\omega \sim N(0, 0.1\sigma^2)$ and $\epsilon=0.02$.
The proportions of selecting the ``true'' models ($E \supseteq \supp{b}$) is 
$0.91$ without randomization and $0.852$ in MAGIC. Much more different
is the power of the two procedures; for a level-$0.05$ test, the statistical 
powers defined above is $0.176$ without randomization and $0.887$ in MAGIC.
Figure \ref{fig:pvalue_alt} is the histograms for the alternative p-values with or
without randomization.

\section{Proofs}
\label{sec:proof}

\subsection{Proof for Theorem \ref{thm:pvalue}}
\begin{proof}
Let $\cS$ be the space for $S$, then $(S, \omega) \in \cS \times \real^p$. 
The joint distribution of $(S, \omega)$ conditional on $(S, \omega) \in \cA$ 
has the following density with respect to the measure $\mu(\diff s)G(\diff \omega)$
\begin{equation}
\label{eq:joint_density}
h(s, \omega) = \frac{\exp\left[b^T T(s)\right]\mathbf{1}\{(s, \omega) \in \cA\}}
{\int_{\cS \times \real^p} \exp\left[b^T T(s)\right] \mathbf{1}\{(s, \omega) \in \cA\} \mu(\diff s)G(\diff\omega)}.
\end{equation}
Since the denominator is merely a normalizing constant, \eqref{eq:joint_density} is also
an exponential distribution with parameters $b$, sufficient statistics $T(S)$ and
a slightly different reference measure $\mathbf{1}\{(s, \omega) \in \cA\}\mu(\diff s)G(\diff \omega)$.
Since $E \supseteq \supp{b}$, $b^T T(s) = b_E^T T_E(s)$, where $b_E, T_E \in \real^{|E|}$
are $b$ and $T$ restricted to set $E$. Thus \eqref{eq:joint_density} can be seen as
an exponential family with sufficient statistics $T_E$ and parameters $b_E$.
To test any hypothesis $H_{0j}: b_j = 0, ~ j \in E$, 
Chapter 4 of \cite{lehman_testing_1997} states that Uniformly Most Powerful Unbiased tests can be
constructed using the statistic $T_j$ and conditioning on all the other sufficient statistics $T_{E\backslash j} \in \real^{|E|-1}$.
Thus the conditional density of the one dimensional distribution for $T_j$ is
\begin{equation}
\label{eq:conditional_density}
\begin{aligned}
h_j(t_j; t_{E\backslash j}) &= \frac{\exp\left[b_j t_j + b_{E\backslash j}^T t_{E\backslash j}\right] \cdot 
\int_{\cA} \mathbf{1}\{T_{E\backslash j}(s) = t_{E\backslash j}\}\mu(\diff s)G(\diff\omega)}
{\exp\left[b_{E\backslash j}^T t_{E\backslash j}\right] \cdot \int_{\cA} \exp\left[b_j T_j(s)\right] 
\mathbf{1}\{T_{E\backslash j}(s) = t_{E\backslash j}\}\mu(\diff s)G(\diff\omega)} \\
&= \frac{\exp(b_j t_j)\int_{\cA} \mathbf{1}\{T_{E\backslash j}(s) = t_{E\backslash j}\}\mu(\diff s)G(\diff\omega)}
{\int_{\cA} \exp\left[b_j T_j(s)\right] 
\mathbf{1}\{T_{E\backslash j}(s) = t_{E\backslash j}\}\mu(\diff s)G(\diff\omega)} 
\end{aligned}
\end{equation}
Thus \eqref{eq:conditional_density} is the density for the distribution 
\begin{equation}
\label{eq:conditional:law}
T_j(S) | T_{E\backslash j}(S), ~(S, \omega) \in \cA,~ (S, \omega) \sim F \times G.
\end{equation}
Note $\eqref{eq:conditional_density}$ involves only the parameter $b_j$, thus it
can be used to test the composite hypothesis $H_{0j}: b_j = 0$, with $b_{E\backslash j}$
taking arbitrary values.
 
Let $H_j$ denote the c.d.f of the above law: $H_j(t_j; T_{E\backslash j}) = \int_{-\infty}^{t_j} h_j(r; T_{E \backslash j}) \diff r$.
Then we can construct our function $\tilde{P}_j: \real^{p} \rightarrow \real$ as 
\begin{equation}
\label{eq:pvalue_construction}
\begin{aligned}
\tilde{P}_j(t) &= \frac{\int_{\cA} \exp(b_j T_j(s))\mathbf{1}\{T_j(s) > t_j\}
\mathbf{1}\{T_{E\backslash j}(s) = t_{E\backslash j}\}\mu(\diff s)G(\diff\omega)}
{\int_{\cA} \exp\left[b_j T_j(s)\right] 
\mathbf{1}\{T_{E\backslash j}(s) = t_{E\backslash j}\}\mu(\diff s)G(\diff\omega)} \\
&= 1-H_j(t_j; t_{E\backslash j}).
\end{aligned}
\end{equation}
Under the null hypothesis, we take $b_j = 0$, thus $\tilde{P}_j$ depends only on $T(s)$ and $\cA$.
We define $P_j(s) = \tilde{P}_j(T(s))$. Now we prove the level-$\alpha$ control (\ref{eq:alpha_level}).
Note
$$
\Pp_{F\times G} \left[P_j(S) \leq \alpha \mid (S, \omega) \in \cA \right]
= \Ee \left[\Pp_{F\times G} [P_j(S) \leq \alpha \mid T_{E\backslash j}, ~ (S, \omega) \in \cA]\right],
$$
it suffices to prove the quantity inside the expectation has the level-$\alpha$ control
for any $T_{E\backslash j}$.
Since $H_j$ is the c.d.f of the conditional law \eqref{eq:conditional:law}, 
$$
\begin{aligned}
&\Pp_{F\times G} [P_j(S) \leq \alpha \mid T_{E\backslash j}, ~ (S, \omega) \in \cA] \\
=&\Pp_{F\times G} [1-H_j(T_j; T_{E\backslash j}) \leq \alpha \mid T_{E\backslash j}, ~ (S, \omega) \in \cA]\\ 
=&\Pp_{F\times G} [T_j \geq H_j^{-1}(1-\alpha) \mid T_{E\backslash j}, ~ (S, \omega) \in \cA]\\ 
=&1-H_j[H_j^{-1}(1-\alpha)] \leq \alpha,
\end{aligned}
$$
where $H_j^{-1}$ generalized inverse for $H_j$, the equality holds when $H_j$ is strictly increasing in $t_j$.
\end{proof}


\subsection{Proof for Lemma \ref{lem:KKT}}
\begin{proof}
Equation (\ref{eq:canonical:random:program}) is a
convex optimization problem. The solution $\hbeta$ and subgradient of 
the $\ell_1$ norm $\hz$ satisfy the Karush–Kuhn–Tucker conditions (KKT),
which are sufficient and necessary.
$$
\begin{cases}
\partial_{\beta}\ell(\hbeta; S) + \hz - \omega = 0, \\
\hz \in \partial \|\hbeta\|_1.
\end{cases}
$$
The equations are simply the differentiation of the optimization objective function.
This gives the equation part in (\ref{eq:KKT}) of the lemma. 
Note the penalty term $\|\hbeta\|_1$ is differentiable except at $0$, its subgradient
at $0$ is $[-1, 1]$. Thus, conditioning on the active set $\hE(S, \omega) = E$ it is
equivalent to: 
$$
\begin{cases}
\hz_j = \text{sign}(\hbeta_j), \quad \forall ~j \in E,\\
|\hz_j| \leq 1, \quad \forall ~j \not\in E.
\end{cases}
$$
Combining the above two, we have the conclusion of the lemma.
\end{proof}

\subsection{Proof for Theorem \ref{thm:change:measure}}
\begin{proof}
Per the discussion above Theorem \ref{thm:change:measure}, it is
not hard to see the distributional constraint on $(S, \hbeta, \hz)$
is that $\gamma(S, \hbeta) + \lambda \hz \sim G$ and is independent of $S$.
Moreover, $(\hbeta, \hz)$ are constrained to be in the region $\cB$.
Thus the law (\ref{eq:change:measure}) is the marginal law of $S$
conditional on selection. 

Now we investigate the joint density of $(s, \beta, z)$.
Through the reconstruction map,
$$
\psi(s, \beta, z) = \bigg(s, \gamma(s, \beta) + \lambda z \bigg),
$$
the density of $(s, \beta, z)$ is simply the product $f(s)g\big(\gamma(s, \beta) + \lambda z \big)$
times the determinant of the Jacobian matrix. 
Standard multivariate calculus yields the form of the Jacobian matrix
of $\psi$ as follows,
$$
\begin{pmatrix}
I & 0 \\
D_s \psi(s,\beta,z) & D_{(\beta,z)} \psi(s,\beta,z) 
\end{pmatrix}
$$
with determinant $\det D_{(\beta,z)} \psi(s,\beta,z)$. Thus we have (\ref{eq:change:measure:density}).

Notice the construction
\eqref{eq:pvalue_construction}, $P_j(s) = 1 - H_j(T_j(s), T_{E\backslash j}(s))$
and $H_j$ is the c.d.f for the conditional distribution \eqref{eq:conditional:law}.
It is equivalent to sampling (\ref{eq:change:measure}) while further conditional on 
$T_{E\backslash j}(S)$. 
After we acquire $m$ samples $\{S^{(1)}, \dots, S^{(m)}\}$, we can approximate the integral
in \eqref{eq:pvalue_construction} as the percentile of $T_j(S_{obs})$ among $\{T_j(S^{(1)}), \dots, T_j(S^{(m)})\}$.
\end{proof}



\subsection{Proof for Corollary \ref{cor:lasso:sampler}}
\begin{proof}
    Notice that once we condition on the active set $E$, and the signs $z_{E,obs}$, 
    \begin{equation}
        \label{eq:reconstruction}
    \beta = \begin{pmatrix}
        \beta_E \\
        0
    \end{pmatrix},
    \qquad 
    z = \begin{pmatrix}
        z_{E, obs} \\
        z_{-E}
    \end{pmatrix}.
    \end{equation}
    Therefore, the density of $(y, \beta, z)$ is equivalent to that of $(y, \beta_E, z_{-E})$, through the construction in \eqref{eq:reconstruction}. 
    Note the Jacobian matrix is 
    $$
    D_{\beta_E, z_{-E}} \psi(y, \beta_E, z_{-E}) =  \begin{pmatrix} X_E^T X_E & 0 \\ X_{-E}^T X_E & \lambda I\end{pmatrix}, 
    $$ 
    where $I$ is the identity matrix of dimension $p-|E|$. 
    Then the Jacobian $J\psi(y,\beta,z) = \lambda^{p-|E|}\det (X_E^T X_E)$. Since the Jacobian 
    is a constant only depending on $X$. 
    Thus, if we plug in $\beta$, $z$ in Theorem \ref{thm:change:measure},
     the density of $(y, \beta_E, z_{-E})$ has the form in Corollary \ref{cor:lasso:sampler}.
\end{proof}

\subsection{Proof for Theorem \ref{thm:logistic}}
\begin{proof}
    We first reconstruct the gradient $\nabla\ell(\hbeta_E)$ from $\hbeta_E$ and 
    $$
    T=
    \begin{pmatrix}
    \bar{\beta}_E \\
     X_{-E}^T( y - \pi(X_E\bar{\beta}_E)) 
    \end{pmatrix}.
    $$
    The Taylor expansion of $\nabla\ell(\hbeta_E)$ at $\bar{\beta}_E$ is 
    $$ 
    \begin{aligned}
    \nabla \ell(\hbeta_E) &= -\frac{1}{\sqrt{n}}X^T(y - \pi(X_E\hbeta_E)) \\
    &= -\frac{1}{\sqrt{n}}X^T\left(y - \pi(X_E\bar{\beta}_E) 
    - W(X_E\bar{\beta}_E)X_E(\hbeta_E-\bar{\beta}_E)\right) + R,
    \end{aligned}
    $$
    where $R = o_p(1)$. Since $\bar{\beta}_E$ is the minimizer of the logistic regression with $E$ variables,
    the gradient at $\bar{\beta}_E$ is zero,
    $$
    X_E^T(y - \pi(X_E\bar{\beta}_E)) = 0. 
    $$
    Thus we can rewrite $\nabla\ell(\hbeta_E)$ in terms of $\hbeta_E$ and $T$ via the following map,
    \begin{equation}
    \label{eq:gradient:recon}
    \nabla \ell(\hbeta_E) = \frac{1}{\sqrt{n}}\begin{pmatrix} X_E^TW(X_E\bar{\beta}_E) X_E(\hbeta_E - \bar{\beta}_E) \\ X_{-E}^TW(X_E\bar{\beta}_E) X_E(\hbeta_E - \bar{\beta}_E) -X_{-E}^T(y-\pi(X_E\bar{\beta}_E)) \end{pmatrix} + R.
    \end{equation}
    Notice that $\bar{\beta}_E$ is the MLE for the negative logistic likelihood, and thus satisfy the asymptotic normality, 
    with asymptotic mean $b_E$, when $E \supseteq \supp{b}$.  Moreover, the following part has asymptotically mean $0$,
    $$
    \begin{aligned}
    &\Ee\left[\frac{1}{\sqrt{n}}X_{-E}^T(y - \pi(X_E\bar{\beta}_E))\right]\\ 
    =& \Ee\left[\frac{1}{\sqrt{n}}X_{-E}^T(y - \pi(X_Eb_E))\right] - \Ee\left[\frac{1}{\sqrt{n}}
    X_E^T W(X_Eb_E) X_E (\bar{\beta}_E - b_E)\right] + o_p(1) \rightarrow 0.
    \end{aligned}
    $$
    Thus we have the asymptotic normality as in (\ref{eq:asymptotic:normality}). Moreover, since
    $\omega = \nabla\ell(\hbeta_E) + \lambda \hz + \epsilon \hbeta_E$,
    then we have asymptotically,
    $$
    (T, \nabla\ell(\hbeta_E) + \lambda \hz + \epsilon \hbeta_E) \overset{d}{\rightarrow} F \times G.
    $$
    The Jacobian is $\det(X_E^T W(X_E\bar{\beta}_E)X_E + \epsilon I)$ which by law of large numbers converges to
    $\det\bigg[\Ee(X_E^T W(X_Eb_E) X_E) + \epsilon I\bigg]$, a constant. Therefore, we have the density (\ref{eq:asymptotic:density}) 
    if we plug in the map \eqref{eq:gradient:recon} for $\nabla\ell(\hbeta)$.
\end{proof}

\subsection{Proof for Corollary \ref{cor:ns_sampling}}
\begin{proof}
    For every node $i$, the $i$-th coordinate of $\beta^i$ is held to be zero, and (\ref{eq:ns_optimize}) is in fact a regression of dimension $p-1$, thus
    $\gamma(X, \hbeta^i) = -X_{-i}^T (x_i - X \hbeta^i) \in \real^{p-1}$,
    and the reconstruction map, 
    $$
    \psi: (X, \hB, \hz) \mapsto (X, \gamma(X, \hB) + \lambda\hz), 
    $$
    where 
    $$
    \begin{aligned}
        &\gamma(X, \hB) = (\gamma(X, \hbeta^1), \dots, \gamma(X, \hbeta^p)) \in \real^{(p-1) \times p}, \\
        &\hz = (\hz^1, \hz^2, \dots, \hz^p),~ \hB = (\hbeta^1, \dots, \hbeta^p),
    \end{aligned}
    $$
    and $\hz^i = \begin{pmatrix} z_{E, obs}^i\\ \hz_{-E}^i\end{pmatrix}$ is the subgradient of the optimization problem (\ref{eq:ns_optimize}). 
    Since $\omega^i$'s are independent, and the Jacobian 
    $$
    J\psi(X,z,B) = \prod_{i \in \Gamma} \det(X_{E^i}^T X_{E^i}),
    $$
    density (\ref{eq:ns_law}) follows.
\end{proof}

\section{Monte-Carlo sampler}
\label{sec:sampling}
Theorem \ref{thm:change:measure} gives an explicit way of computing the density for the law of selective inference. We can 
use a Gibbs sampler to rotate through sampling $(S, \hbeta, \hz)$. For sampling $S$ and $\hbeta$, we can take a Metropolis-Hastings
step and use the density to compute the acceptance probability. For sampling $\hz$, it is even simpler as we recognize the conditional
distribution of $\hz|S, \hbeta$ is simply a truncated $G$ distribution. To illustrate our sampler, we describe the sampling scheme of 
some of our examples in more details.

\subsection{Randomized Lasso sampler}
\label{sec:lasso:sampler}
Without loss of generality, we assume the density of added noise $G$ is symmetric and each coordinate of $\omega$ is independent.
This is in fact what we use a lot in practice. 
Also denote $G_{\Delta^-, \Delta^+}$ as truncated distribution $G$ with $\Delta^-$, $\Delta^+$ as the lower and upper truncation points,
and $h(y, \hbeta_E, \hz_{-E})$ to be the density in (\ref{eq:lasso:density}). 
Then to test the null hypothesis $H_{0j}: b_j = 0$, we propose Algorithm \ref{alg:lasso}. Note the step sizes $a_n$ and $c_n$ in
Algorithm \ref{alg:lasso} is chosen through \cite{robust_metropolis} to achieve the desired acceptance rate. 

\begin{algorithm}[tb]
    \caption{Metropolis Hastings sampler for randomized Lasso}
    \label{alg:lasso}
\begin{algorithmic}
    \STATE {\bf Set:} $b=0$ for distribution $f_b$, compute the explicit expression $h$.
    \STATE {\bf Compute:} $P = X_{E\backslash j} X_{E\backslash j}^{\dagger}$, $R = I - P$,
    \STATE {\bf Initialize:} $(y^0, \hbeta_E^0, \hz_{-E}^0) \leftarrow (y, \hbeta_E, \hz_{-E})$,
    \STATE {\bf Step data:} $y^{(n+1)} \leftarrow Py^{(n)} + a_n \cdot R\tau, \quad \tau \sim N(0, I)$,
    compute the acceptance ratio $r = \frac{h(y^{(n+1)}, \hbeta_E^{(n)}, \hz_{-E}^{(n)})}{h(y^{(n)}, \hbeta_E^{(n)}, \hz_{-E}^{(n)})}$,
    accept $y^{(n+1)}$ with probability $r$, otherwise $y^{(n+1)} \leftarrow y^{(n)}$. If $r > 1$, accept $y^{(n+1)}$. 
    \STATE {\bf Step coefficient:} $\hbeta_E^{(n+1)} \leftarrow s_E |\hbeta_E^{(n)} + c_n \cdot \nu|, ~ \nu \sim G$,
    compute the acceptance ratio $r = \frac{g(y^{(n+1)}, \hbeta_E^{(n+1)}, \hz_{-E}^{(n)})}{g(y^{(n+1)}, \hbeta_E^{(n)}, \hz_{-E}^{(n)})}$,
    and accept/reject accordingly.
    \STATE {\bf Step subgradient:} compute the upper and lower limits,
    $$
    \begin{aligned}
        \Delta^+ &= -X_{-E}^T(y^{(n+1)} - X_E \hbeta_E^{(n+1)}) + \lambda \mathbf{1}, \\ 
        \Delta^- &= -X_{-E}^T(y^{(n+1)} - X_E \hbeta_E^{(n+1)}) - \lambda \mathbf{1}, 
    \end{aligned}
    $$
    sample $\lambda \hz_{-E}^{(n+1)} \overset{ind}{\sim} G_{\Delta^-, \Delta^+}$.
\end{algorithmic}
\end{algorithm}

\subsection{Neighborhood selection}
Similar to the scheme in Section \ref{sec:lasso:sampler}, we use a Gibbs sampler to sample $X$, $\hB$ and $\hz$ respectively. The
sampling for $\hB$ and $\hz$ are analogous to that of Section \ref{sec:lasso:sampler}, 
and we only need a proposal distribution for $X$.
As mentioned in Section \ref{sec:ns}, to test the hypothesis $H_{0, ij}: \Theta_{ij} = 0$, we condition on $\{x_{i'}^T x_{j'}, (i',j') \neq (i,j)\}$. 
To sample the data matrix $X$, we rotate through its columns, sampling one column at a time, keeping all the others as constant.
More specifically, for column $i$, we sample from the distribution,
$$
x_i | X_{-i}, \|x_i\|^2, ~ x_{i'}^T x_{j'}, ~(i', j') \in E^{\lor}, (i', j') \neq (i,j).
$$
Note the graph structure gives a natural partition of the nodes into different connected components, let $\text{ne}(i)$ be
the nodes in the connected component of $i$, then 
$x_i \perp x_k, ~\forall k \not\in \text{ne}(i)$, conditioning on all the other $x_j$' in $\text{ne}(i)$.
Thus the above law is equivalent to,
\begin{equation}
    \label{eq:ns:column:law}
    x_i | \|x_i\|^2, x_j, x_i^T x_j, j \in \text{ne}(i).
\end{equation}
We can sample the above law \eqref{eq:ns:column:law} by sampling uniformly from a sphere with radius $\|x_i\|$, holding 
the projections onto the $x_j$'s constant. 
After sampling a new column of X, we compute the accept ratio, accept/reject accordingly 
and move to the next column. As for the sampling of $\hB$ and $\hz$, we can develop an algorithm similar 
to Algorithm \ref{alg:lasso}.

\section{Discussion}
\label{sec:discussion}

MAGIC has the following limitations that we hope to remove in future work. First, the penalty
in our convex program have to be $\ell_1$ penalty. Second, we assume parametric models, more
specifically in the exponential family setting. Third, in the setting for Section \ref{sec:logistic_lasso},
we require the dimension $p$ to be fixed, leaving the high-dimensional problem $p>n$ as an interesting problem.  

\bibliographystyle{plain}
\bibliography{sampling}
\end{document}